\documentclass[lettersize,journal]{IEEEtran}
\usepackage{amsmath,amsfonts}
\usepackage{algorithmic}
\usepackage{algorithm}
\usepackage{array}
\usepackage[caption=false,font=normalsize,labelfont=sf,textfont=sf]{subfig}
\usepackage{textcomp}
\usepackage{stfloats}
\usepackage{url}
\usepackage{verbatim}
\usepackage{graphicx}
\usepackage{subcaption}
\usepackage{cite}
\usepackage{booktabs}
\usepackage{multirow}
\usepackage{hyperref}
\usepackage{csquotes}
\usepackage{color}

\usepackage{tabularx}
\usepackage{amsmath}
\usepackage[intoc]{nomencl}
\makenomenclature
\usepackage{etoolbox}
\renewcommand\nomgroup[1]{%
  \item[\bfseries
  \ifstrequal{#1}{A}{Indices}{%
  \ifstrequal{#1}{B}{Sets}{
  \ifstrequal{#1}{C}{Parameters}{%
  \ifstrequal{#1}{D}{Variables}{}}}}%
]}
\usepackage{makecell}
\usepackage{colortbl}
\usepackage[table]{xcolor} 
\usepackage{mathtools}
\definecolor{highlightcol}{rgb}{0.85, 0.9, 1.0} 
\begin{document}


\title{Relax-and-Cut for Temporal SCUC Decomposition}

\author{
Jinxin Xiong,~\IEEEmembership{Student Member,~IEEE,}
Linxin Yang,~\IEEEmembership{Student Member,~IEEE,}
Yingxiao Wang,~\IEEEmembership{Student Member,~IEEE,}
Yanting Huang,~\IEEEmembership{Student Member,~IEEE,}
Jianghua Wu,~\IEEEmembership{Member,~IEEE,}
Shunbo Lei,~\IEEEmembership{Senior Member,~IEEE,}
Akang Wang,~\IEEEmembership{Member,~IEEE}

\thanks{This work is supported by Natural Science Foundation of China (Grant No.
12301416), Guangdong Basic and Applied Basic Research Foundation (Grant No. 2024A1515010306), Shenzhen Science and Technology Program (Grant No. RCBS20221008093309021), and Hetao Shenzhen-Hong Kong Science and Technology Innovation Cooperation Zone Project (No. HZQSWS-KCCYB-2024016). \textit{(Corresponding author: Akang Wang)}}
\thanks{Jinxin Xiong, Linxin Yang and Akang Wang are with School of Data Science, The Chinese University of Hong Kong, Shenzhen; Shenzhen International Center for Industrial and Applied Mathematics, Shenzhen Research Institute of Big Data, China (e-mail: jinxinxiong, linxinyang@link.cuhk.edu.cn; wangakang@sribd.cn).}
\thanks{Yingxiao Wang, Yanting Huang and Shunbo Lei are with School of Science and Engineering, The Chinese University of Hong Kong, Shenzhen; Shenzhen Research Institute of Big Data, China (e-mail: wangbottlecap@gmail.com; yantinghuang@link.cuhk.edu.cn; leishunbo@cuhk.edu.cn).}
\thanks{Jianghua Wu is with Shenzhen Research Institute of Big Data, China (e-mail: wujianghua@sribd.cn),}
}


\maketitle
\begin{abstract}
The Security-Constrained Unit Commitment (SCUC) problem presents formidable computational challenges due to its combinatorial complexity, large-scale network dimensions, and numerous security constraints. 
While conventional temporal decomposition methods achieve computational tractability through fixed short-term time windows, this limited look-ahead capability often results in suboptimal, myopic solutions.
We propose an innovative \emph{relax-and-cut} framework that alleviates these limitations through two key innovations. First, our enhanced temporal decomposition strategy maintains integer variables for immediate unit commitment decisions while relaxing integrality constraints for future time periods, thereby extending the optimization horizon without compromising tractability. Second, we develop a dynamic cutting-plane mechanism that selectively incorporates N-1 contingency constraints during the branch-and-cut process, avoiding the computational burden of complete upfront enumeration. The framework optionally employs a Relaxation-Induced Neighborhood Search procedure for additional solution refinement when computational resources permit.
Comprehensive numerical experiments demonstrate the effectiveness of our approach on large-scale systems up to 13,000 buses. The proposed method can achieve optimality gaps below 1\% while requiring only 20\% of the computation time of monolithic Gurobi solutions. Compared to existing decomposition approaches, our framework provides superior performance, simultaneously reducing primal gaps by 60\% and doubling solution speed. These significant improvements make our method particularly well-suited for practical SCUC implementations where both solution quality and computational efficiency are crucial.
\end{abstract}


\nomenclature[A]{$g$}{Generator index.}
\nomenclature[A]{$t$}{Time period index.}
\nomenclature[A]{$b$}{Bus index.}
\nomenclature[A]{$l$}{Transmission line index.}
\nomenclature[A]{$c$}{Contingency index, $0$ for the base case.}

\nomenclature[B]{$\mathcal{G}$}{Set of generators.}
\nomenclature[B]{$\mathcal{B}$}{Set of buses.}
\nomenclature[B]{$\mathcal{G}_b$}{Set of generators connecting to bus $b$.}
\nomenclature[B]{$\mathcal{T}$}{Set of time periods $1,\cdots, T$.}
\nomenclature[B]{$\mathcal{L}$}{Set of transmission lines.}
\nomenclature[B]{$\mathcal{C}$}{Set of contingencies.}
\nomenclature[B]{$[m,n]$}{Set of integers $m, m+1, \cdots, n$.}

\nomenclature[C]{$\bar{x}_{g0}$}{Initial status of generator $g$.}
\nomenclature[C]{$\overline{UT}_g$}{Initial minimum up time of generator $g$.}
\nomenclature[C]{$\overline{DT}_g$}{Initial minimum down time of generator $g$.}
\nomenclature[C]{$\overline{p}_{g0}$}{Initial power generation of generator $g$.}
\nomenclature[C]{$UT_g$}{Minimum up time of generator $g$.}
\nomenclature[C]{$DT_g$}{Minimum down time of generator $g$.}
\nomenclature[C]{$P^L_g$}{Minimum power limit of generator $g$.}
\nomenclature[C]{$P^U_g$}{Maximum capacity limit of generator $g$.}
\nomenclature[C]{$SU_g$}{Start-up capacity of generator $g$.}
\nomenclature[C]{$SD_g$}{Shut-down capacity limit of generator $g$.}
\nomenclature[C]{$RU_g$}{Maximum ramp up capacity of generator $g$.}
\nomenclature[C]{$RD_g$}{Maximum ramp down capacity of generator $g$.}
\nomenclature[C]{$D_{bt}$}{Demand of bus $b$ in time period $t$.}
\nomenclature[C]{$F_{l}^c$}{Thermal limit of line $l$ under 
contingency $c$.}
\nomenclature[C]{$C^{\text{curtail}}$}{Penalty related to load curtail.}
\nomenclature[C]{$C^{SU}_g$}{Startup cost for generator $g$.}
\nomenclature[C]{$C^{NL}_g$}{No-load cost for generator $g$.}

\nomenclature[D]{$x_{gt}$}{Commitment status of generator $g$ in time period $t$.}
\nomenclature[D]{$z_{gt}$}{Start-up status of generator $g$ in time period $t$.}
\nomenclature[D]{$w_{gt}$}{Shut-down status of generator $g$ in time period $t$.}
\nomenclature[D]{$p_{gt}$}{Power output of generator $g$ in time period $t$.}
\nomenclature[D]{$p_{gt}$}{Power output of generator $g$ in time period $t$.}
\nomenclature[D]{$\delta_{lb}^c$}{Power Transfer Distribution Factor (PTDF) of bus $b$ to line $l$ under contingency $c$.}
\nomenclature[D]{$d_{bt}^{\text{curtail}}$}{Amount of load curtailed for bus $b$ at time $t$.}
\nomenclature[D]{$d_{bt}^{\text{param}}$}{Variable for load in the parameterized model for bus $b$ at time $t$.}
\printnomenclature

\section{Introduction}

\IEEEPARstart{S}{ecurity}-Constrained Unit Commitment (SCUC) is a fundamental optimization problem in power system operations, tasked with determining the most cost-effective schedule of generator commitments (on/off states) and dispatch levels while ensuring reliable power delivery over a multi-period horizon. 
In practice, this complex decision-making process must be completed routinely for large-scale systems within tight operational timeframes.

Mathematically, the SCUC problem is formulated as a large-scale \textit{mixed-integer linear program} (MILP). Its computational complexity stems from two primary challenges: (i) a high-dimensional decision space, and (ii) numerous dense security constraints.
The interplay between these factors renders even the initial \textit{linear programming}~(LP) relaxation computationally intensive, failing to meet practical solution time requirements. Consequently, developing advanced computational techniques for rapid identification of high-quality SCUC solutions is critically important.

Significant progress in SCUC solution methods has been achieved through advances in general-purpose MILP solvers and optimization techniques. Modern approaches employ tight MILP formulations~\cite{panConvexHullsUnit2017, knueven2020mixed, pan2016polyhedral, tejada2019unit} that, when integrated with advanced branch-and-cut frameworks, have substantially expanded the scale of tractable SCUC problems. These methods are further enhanced through sophisticated primal heuristics and specialized cutting planes.
Nevertheless, solving large-scale SCUC problems within strict operational time windows remains fundamentally challenging due to the sheer volume of security constraints and long planning horizons requiring temporal coupling.

The explosion of temporal complexity has spurred the development of decomposition-based strategies.
The works of~\cite{safdarian2019temporal,kim2018temporal} proposed \textit{temporal decomposition} strategies that divide the scheduling horizon into subhorizons, formulating a distributed SCUC framework with coupling intervals or auxiliary variables to manage inter-temporal constraints. 
In addition to temporal decomposition, it is also possible to reduce the computational complexity of SCUC by spatially decoupling the problem into geographically or structurally independent subproblems. 
These methods partition the power system network into regions or zones and iteratively optimize local commitment and dispatch decisions while enforcing system-wide coordination to satisfy network security constraints.
For instance, \cite{xavier2024decomposable} proposed a novel formulation for decentralized optimization that ensures efficient management of intra-zonal resources and boundary coordination between subsystems, significantly enhancing the scalability and structural sparsity of large-scale systems.
However, these decentralized methods typically require coordination among subproblems, which can potentially lead to slow convergence or necessitate careful tuning of auxiliary parameters to ensure solution quality.

To mitigate the explosion of security constraints, prior studies~\cite{zhai2010fast,ardakani2014acceleration}
proposed methods to identify critical constraints.
However, these methods faced scalability limitations due to the need to solve auxiliary optimization problems.
Subsequent works~\cite{chen2016improving, tejada2017security} improved efficiency through iterative security constraint filtering, starting with a relaxed problem that excludes the vast majority of transmission constraints and incrementally adding violated constraints until no violations remain.
~\cite{bouffard2005umbrella, ardakani2013identification} observed that enforcing a small subset of security constraints is usually sufficient to ensure the satisfaction of all the remaining constraints.
Leveraging this insight,~\cite{xavier2019transmission} enhanced scalability by selectively adding only the most violated constraints in each iteration, substantially reducing problem complexity. 
Although these approaches incrementally reduce problem size, a significant drawback is the need to solve the MILP model from scratch in every iteration, delaying the availability of a feasible solution until the entire process concludes.

\begin{figure*}[htbp]
    \centering
\includegraphics[width=0.95\linewidth]{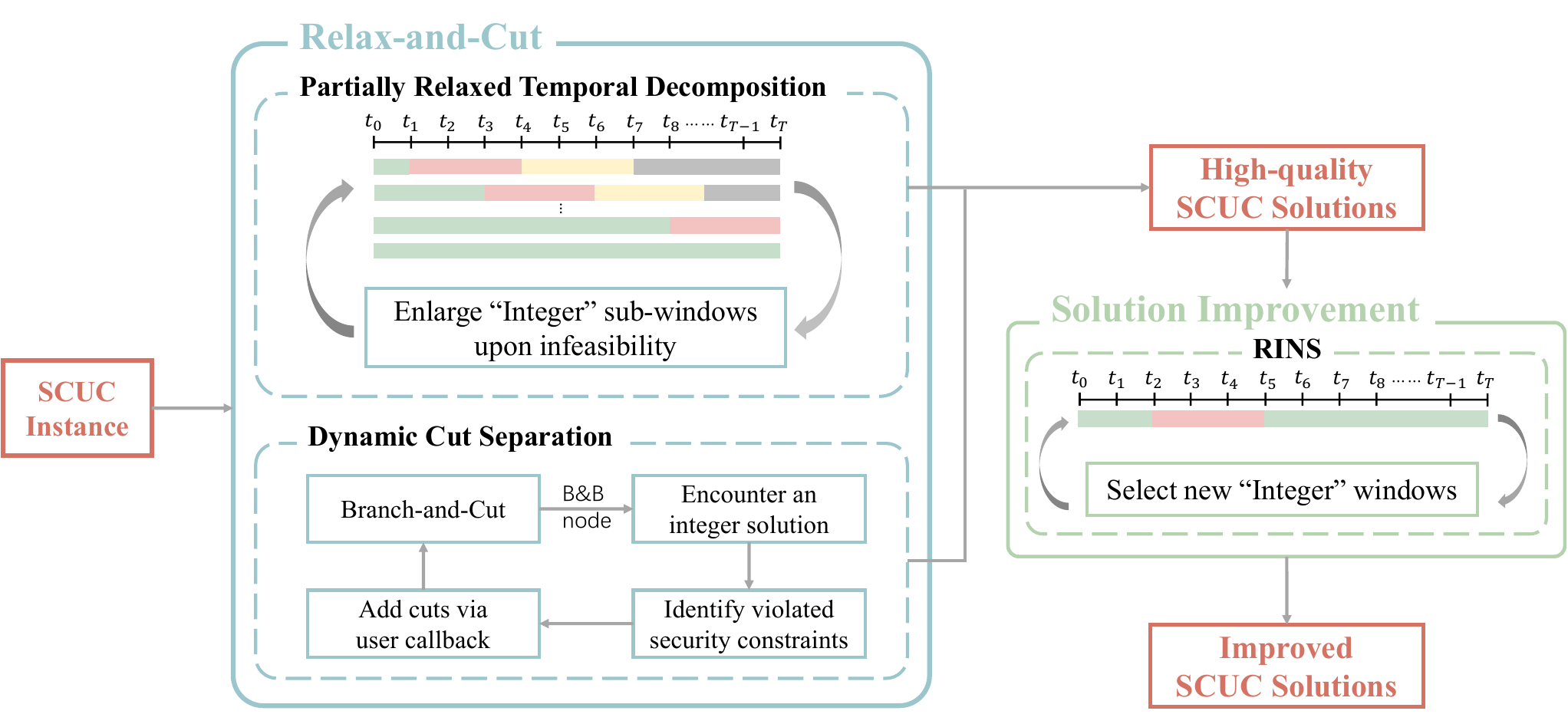}
    \caption{
    A flowchart of the proposed relax-and-cut framework. The process begins with a partially relaxed temporal decomposition to rapidly find a high-quality feasible solution, which can then be further refined using RINS. Throughout the framework, security constraints are dynamically enforced within a branch-and-cut procedure.}
    \label{fig:flowchart}
\end{figure*}

This paper introduces a novel \textit{relax-and-cut} method to enhance temporal decomposition for efficiently identifying high-quality SCUC solutions. 
Specifically, the proposed framework advances through the planning horizon by solving a sequence of subproblems. 
In each subproblem, commitment decisions are managed through a unique partitioning scheme—categorizing discrete variables as fixed from prior steps, constrained in the current window, relaxed in a look-ahead period, or deferred to future subproblems. This strategy yields a series of computationally tractable yet informative, partially relaxed subproblems. To handle the combinatorial challenge of system security, N-1 contingency constraints are not enumerated upfront but enforced dynamically as cutting planes within the native \textit{branch-and-bound}~(B\&B) process.
The integrated relax-and-cut framework, combining the partially relaxed temporal decomposition strategy and dynamic cut separation, alleviates the fundamental trade-off between solution quality and computational complexity.
Upon achieving an initial feasible solution from the relax-and-cut framework, the Relaxation-Induced Neighborhood Search (RINS)~\cite{danna2005exploring} is adapted to systematically explore the solution's vicinity for further improvements. This integrated methodology, illustrated in Fig.~\ref{fig:flowchart}, provides a robust and scalable approach to large-scale SCUC problems.

The distinct contributions of this paper are summarized as follows:
\begin{itemize}
\item \textbf{Partially Relaxed Temporal Decomposition.} We enhance temporal decomposition through strategic partial relaxation of integrality constraints, selectively applying continuous relaxations to unit commitment decisions in future time windows while maintaining exact MILP formulations for immediate commitments. 

\item \textbf{Dynamic Cut Separation.} 
We propose to integrate security constraint enforcement directly within the B\&B process, combining parallel constraint filtering with solver multi-threading for synergistic performance gains.

\item \textbf{Superior Performances.} 
Through extensive computational experiments on public benchmarks—including large-scale instances with over 13,000 buses—we demonstrate that our method achieves sub-1\% optimality gaps in merely 20\% of the runtime required by monolithic SCUC solutions using Gurobi, while simultaneously reducing primal gaps by 60\% and halving solution times compared to existing temporal decomposition approaches.

\end{itemize}

The paper is organized as follows: Section~\ref{sec:problem_def} formulates the SCUC problem, followed by our partially relaxed temporal decomposition in Section~\ref{sec:framework}. Section~\ref{sec:dynamic_cut} presents dynamic enforcement for security constraints, while Section~\ref{sec:overall_algo} integrates these into the relax-and-cut framework with RINS refinement. Section~\ref{sec:case_study} evaluates performance against state-of-the-art methods. We conclude with a discussion and future directions in Section~\ref{sec:conclusion}. 

For reproducibility, our implementation is publicly available at~\url{https://github.com/jx-xiong/relax-and-cut.git}.

\section{Problem Definition}
\label{sec:problem_def}
The SCUC problem determines the most economical generation schedule while satisfying operational constraints over a finite time horizon $\mathcal{T} \coloneqq \left\{1, 2, ..., T\right\}$. 
To elucidate the fundamental computational challenges of SCUC, this section introduces a simplified formulation that retains only the essential problem components. 
For complete mathematical modeling details and constraint formulations, we refer interested readers to~\cite{knueven2020mixed}.

\noindent\textit{Initial status:}
\begin{equation}
\label{eq:initial_status}
    \begin{aligned}
    & x_{gt} = 1 \quad \forall g\in \mathcal{G}: \bar{x}_{g0}=1, t\in [0, \overline{UT}_g]\\
    & x_{gt} = 0 \quad \forall g\in\mathcal{G}: \bar{x}_{g0}=0, t\in [0, \overline{DT}_g] 
\end{aligned}
\end{equation}

\noindent\textit{Logic constraints:}
\begin{equation}
\label{eq:status}
\begin{aligned}
    & x_{gt}-x_{g,t-1}=z_{gt}-w_{gt} && \forall g\in\mathcal{G}, t\in \mathcal{T} \\
    & z_{gt} + w_{gt} \leq 1 && \forall g\in\mathcal{G}, t\in \mathcal{T}  \\
\end{aligned}
\end{equation}

\noindent\textit{Minimum up/down:}
\begin{equation}
\label{eq:minup/down}
    \begin{aligned}
    & \sum_{\tau = \max{\{1, t - UT_g+1 \}}}^{t} z_{g\tau} \leq x_{gt} && \forall g \in \mathcal{G}, t \in \mathcal{T} \\
    & \sum_{\tau= \max{\{1, t - DT_g + 1\}}}^{t} w_{g\tau} \leq 1 - x_{gt} && \forall g \in \mathcal{G}, t \in \mathcal{T} 
\end{aligned}
\end{equation}

\noindent\textit{Initial production:}
\begin{equation}
    \label{eq:init_power}
     p_{g0} = \overline{p}_{g0} \;\; \forall g \in \mathcal{G} 
\end{equation}

\noindent\textit{Production limits:}
\begin{equation}
    \label{eq:power_limit}
    P_{g}^Lx_{gt} \leq p_{gt} \leq P_{g}^Ux_{gt} \;\; \forall g\in\mathcal{G}, t\in\mathcal{T}
\end{equation}

\noindent\textit{Ramp up/down:}
\begin{equation}
    \label{eq:ramp up/down}
    \begin{aligned}
    & p_{gt} - p_{g,t-1} \leq RU_gx_{g,t-1} + SU_gz_{gt}  && \forall g\in\mathcal{G}, t\in\mathcal{T} \\
    & p_{g,t-1} - p_{gt} \leq RD_gx_{gt} + SD_g w_{gt}  && \forall g\in\mathcal{G}, t\in\mathcal{T}
    \end{aligned}
\end{equation}

\noindent\textit{System-wide constraints:}
\begin{equation}
    \begin{aligned}
    & \sum_{g\in\mathcal{G}} p_{gt} = \sum_{b\in\mathcal{B}} \left(d_{bt}^{\text{param}} - d_{bt}^{\text{curtail}} \right) \quad\quad\quad \forall t \in \mathcal{T}\\
    & -F_l^c \leq \sum_{b\in\mathcal{B}} \delta_{lb}^c \left(\sum_{g\in\mathcal{G}_b} p_{gt} - \left(d_{bt}^{\text{param}}-d_{bt}^{\text{curtail}}\right)\right) \leq F_l^c\\
    & \hspace{4cm}\forall l\in\mathcal{L}, t\in\mathcal{T}, c\in\{0\}\cup\mathcal{C}  \\
    & d_{bt}^{\text{param}} = D_{bt} \quad \forall b\in\mathcal{B}, t \in \mathcal{T}
\end{aligned}
\label{eq:security}
\end{equation}

By defining $\Theta \coloneq \left\{\bar{x}_{g0}, \overline{UT}_g, \overline{DT}_g, \bar{p}_{g0}\right\}_{g \in \mathcal{G}}$, the SCUC problem can be formulated as follows:
\begin{equation}
\label{eq:param_prob}
    \begin{aligned}
\text{SCUC}&(T, \left\{D_{bt}\right\}_{b\in\mathcal{B},t \in \mathcal{T}},\Theta ) \\ 
\coloneq &\min \quad \sum_{g\in\mathcal{G}}\sum_{t\in\mathcal{T}}\left(C^P(p_{gt}) + C_g^{SU}z_{gt} + C_g^{NL}x_{gt}\right) \\
&+ C^{\text{curtail}}\sum_{b\in \mathcal{B}} \sum_{t\in\mathcal{T}} d_{bt}^{\text{curtail}} \\
 &\text{s.t.} \;\; \text{constraints~(\ref{eq:initial_status})~--~(\ref{eq:security})} \\
& \;\;\;\quad x_{gt}, z_{gt}, w_{gt} \in \{0,1\} \quad \quad \forall g\in\mathcal{G}, t\in\mathcal{T}.
\end{aligned}
\end{equation}

As shown in~(\ref{eq:param_prob}), the SCUC problem aims to minimize total system costs, including production costs represented by the piecewise linear function $C^P(p_{gt})$ and penalties for demand curtailment. 
The production cost function $C^P(p_{gt})$ admits multiple mathematical representations, as rigorously analyzed in~\cite{knueven2020mixed}. 
The SCUC calls for minimization of the objective function subject to the full set of operational constraints (\ref{eq:initial_status})~--~(\ref{eq:security}), which covers unit commitment rules, economic dispatch requirements, and system security conditions.

\section{Partially Relaxed Temporal decomposition}
\label{sec:framework}

The SCUC formulation in Section~\ref{sec:problem_def}, while fundamental, highlights a critical scalability challenge. 
The number of variables and constraints in (\ref{eq:param_prob}) grows combinatorially with the system size—spanning generators, transmission lines, contingencies, and time periods—rendering monolithic solution approaches computationally challenging for large-scale networks. 
Although conventional temporal decomposition approaches, such as~\cite{barrows2014time}, aim to alleviate this burden, they face a fundamental trade-off: attempts to improve solution quality by extending the optimization window~\cite{kim2018temporal} or creating overlapping subproblems dramatically increase computational complexity, diminishing the very advantage of decomposition.

\begin{figure}[htbp]
    \centering
    \includegraphics[width=1\linewidth]{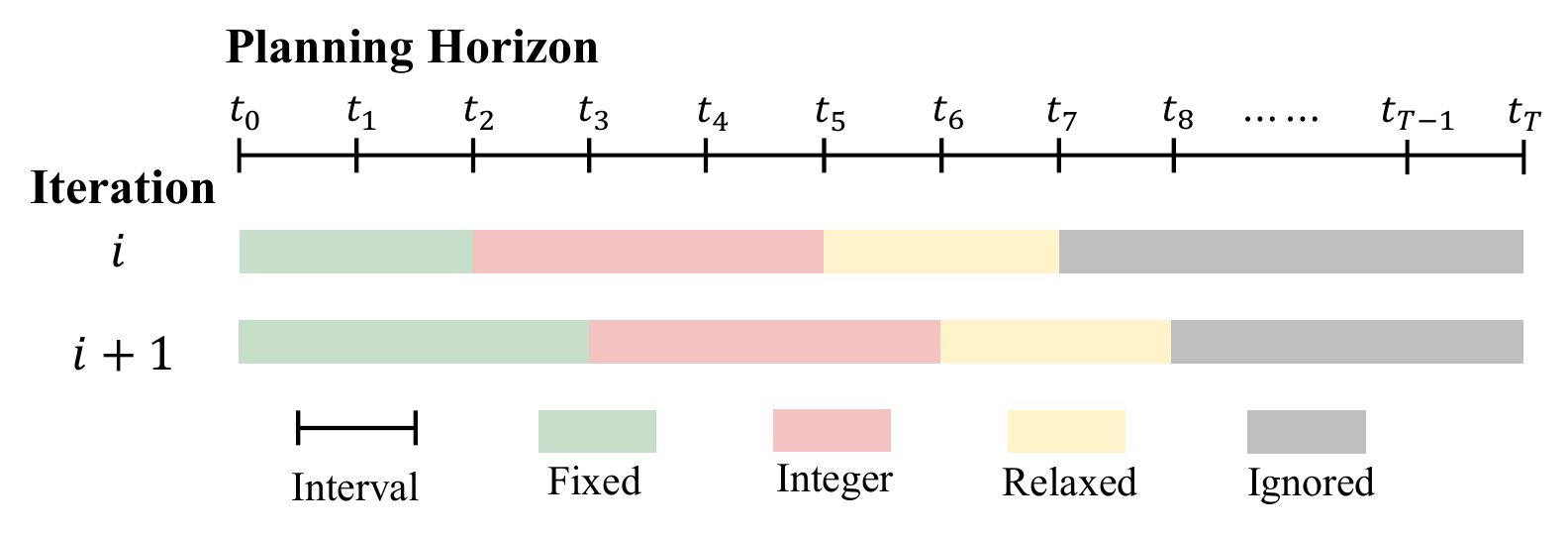}
    \caption{An illustration of the partially relaxed temporal decomposition approach.
    For notational convenience, we denote the entire planning horizon as $[t_0, t_T]$.
    In this example, $s_I=3, s_R=2, \Delta t = 1$.
    At the $i^{th}$ iteration, $w_F$ spans intervals $[t_0, t_2]$, $w_I$ covers $(t_2, t_5]$, and $w_R$ includes $(t_5, t_7]$.
    After solving the $i^{th}$ iteration, the window advanced by $\Delta t$. 
    The interval $(t_2, t_3]$ from previous $w_I$ is merged into $w_F$.
    Also, $w_I$ shifts to $(t_3, t_6]$ and $w_R$ advances to $(t_6, t_8]$.
    }
    \label{fig:rolling}
\end{figure}

However, careful examination of SCUC's temporal characteristics reveals that unit commitment decisions demonstrate only weak coupling across sufficiently separated time periods (e.g., between the first time interval and the tenth). 
This structural insight fundamentally motivates our partially relaxed temporal decomposition strategy. 
Inspired by recent advances~\cite{gasse2022machine,wang2022efficient}, our method strategically limits the optimization window to balance computational tractability with solution quality.
The proposed strategy employs a dynamic temporal decomposition scheme that progressively solves the original problem through advancing time windows (see Fig.~\ref{fig:rolling} for an example). At each step, the planning horizon is partitioned into four distinct temporal regions, where commitment decisions are treated as: (i) \textit{fixed} (for previously determined intervals), (ii) \textit{integer-constrained} (actively optimized decisions), (iii) \textit{relaxed} (continuous relaxation), or (iv) \textit{deferred} (excluded from current optimization). This classification depends on each variable's temporal proximity to the current decision window.

The \enquote{Fixed} window $w_F$ comprises $s_F$ time intervals for which commitment decisions retain their values from prior iterations. Initially, $w_F = \emptyset$.
The \enquote{Integer} window $w_I$ spans $s_I$ time periods where decision variables must satisfy all relevant constraints, encouraging the subproblem to contribute to a feasible SCUC solution.
Discrete variables in the \enquote{Relaxed} window $w_R$ (length $s_R$) are relaxed to continuous values but retain temporal coupling constraints. 
The remaining intervals form the \enquote{Deferred} window, where both variables and constraints are excluded to decouple distant time periods. 
At each iteration, the solved subproblem only entails \enquote{Integer} and \enquote{Relaxed} windows, balancing immediate feasibility with longer-term planning consideration.
The inclusion of $w_R$ mitigates myopic decision-making by extending the solution's lookahead horizon beyond the integer-constrained window $w_I$, thereby improving solution quality while maintaining tractability.

After solving the current subproblem, the window advances by $\Delta t$ intervals.
The \enquote{Fixed} window expands by incorporating $\Delta t$ intervals from the preceding \enquote{Integer} sub-window, while the \enquote{Integer} and \enquote{Relaxed} windows shift forward by $\Delta t$ intervals.
The parameters—such as the initial status and generation level—are properly updated, forming a new subproblem for the next iteration.

By utilizing the parameterized problem~(\ref{eq:param_prob}), we define the partially relaxed subproblem, $\text{SCUC}^{PR}$, as in~(\ref{eq:param_subprob}).
The primary difference between the original SCUC problem in~(\ref{eq:param_prob}) and the $\text{SCUC}^{PR}$ subproblem in~(\ref{eq:param_subprob}) are: (i) \textit{Reduced Time Horizon}: Each subproblem spans a shorter interval, $[1,s_I+s_R]$, which is usually significantly smaller than the entire planning horizon, $[1,T]$;
(ii) \textit{Partial Relaxation}: The discrete decision variables $(x_{g\cdot}, z_{g\cdot}, w_{g\cdot})$ are relaxed to be continuous for the latter part of the subproblem's horizon, $[s_I+1,s_I+s_R]$. This relaxation substantially reduces the problem's combinatorial complexity.
\begin{equation}
\resizebox{0.9\linewidth}{!}{$\displaystyle
\label{eq:param_subprob}
    \begin{aligned}
\text{SCUC}^{PR}&\left(s_I+s_R, \left\{D_{bt}\right\}_{b\in\mathcal{B},t \in w_I\cup w_R}, \Theta\right) \\ 
\coloneq \min &\quad \sum_{g\in\mathcal{G}}\sum_{t\in [1, s_I+s_R]}\left(C^P(p_{gt})+ C_g^{SU}z_{gt} + C_g^{NL}x_{gt}\right) \\
& + c^{\text{curtail}}\sum_{b\in \mathcal{B}} \sum_{t\in [1,s_I+s_R]} d_{bt}^{\text{curtail}} \\
 \text{s.t.} \;\; &\text{constraints~(\ref{eq:initial_status})~--~(\ref{eq:security})} \\
& x_{gt}, z_{gt}, w_{gt} \in \{0,1\}\;\; \forall g\in\mathcal{G}, t\in [1,s_I] \\
& x_{gt}, z_{gt}, w_{gt} \in [0,1]\;\; \forall g\in\mathcal{G}, t\in [s_I+1, s_I+s_R].
\end{aligned}
$}
\end{equation}

The decomposition of the planning horizon into sequential subproblems preserves their inherent coupling through physical constraints, particularly minimum up/down-time requirements~(\ref{eq:minup/down}) and ramp-rate limits~(\ref{eq:ramp up/down}), which create essential temporal dependencies between adjacent subproblems.
Our methodology maintains these critical linkages through systematic state propagation. Each subproblem $i$ is initialized with parameters $\Theta^i \coloneqq \{\bar{x}^i_{g0}, \overline{UT}^i_g, \overline{DT}^i_g, \bar{p}^i_{g0}\}_{g\in\mathcal{G}}$, where $\bar{x}^i_{g0}$ represents the initial commitment status, $\overline{UT}^i_g$ and $\overline{DT}^i_g$ denote remaining required uptime and downtime durations respectively, and $\bar{p}^i_{g0}$ specifies the initial power output for each generator $g\in\mathcal{G}$.
The subproblem solution yields two key outputs: commitment decisions $\{\bar{x}^i_{gt}\}_{g\in\mathcal{G}}$ and generation levels $\{\bar{p}^i_{gt}\}_{g\in\mathcal{G}}$. From these outputs, we derive the actual achieved operational durations ($\widehat{UT}^i_g$ for uptime and $\widehat{DT}^i_g$ for downtime) along with terminal generation states at time $\Delta t$. These computed values serve as inputs for subsequent subproblems.
The state transition mechanism (Algorithm~\ref{alg:update}) ensures both temporal consistency across subproblem boundaries and physical feasibility through strict enforcement of minimum runtime constraints. This partially relaxed temporal decomposition framework successfully achieves dual objectives: it retains the mathematical fidelity of the original monolithic problem while enabling computationally tractable decomposition.
We remark that the framework will adaptively reinitialize with an enlarged \enquote{Integer} window once a subproblem is declared as infeasible.
The inclusion of the \enquote{Relaxed} window mitigates greedy decision-making at each time step, thereby improving feasibility.

\begin{algorithm}[h]
\caption{Initialization for generator $g$}\label{alg:update}
\begin{algorithmic}[1]
\STATE \textbf{Input:} $\Delta t$, $UT_g$, $DT_g$; commitment decisions $\bar{x}_{gt}^i$, generation levels $\bar{p}_{gt}^i$ and consecutive up/down duration $\widehat{UT}_{g}^i / \widehat{DT}_{g}^i$.
\STATE \textbf{Initialize:} $\overline{UT}_g^{i+1} \gets 0$, $\overline{DT}_g^{i+1} \gets 0$.
\STATE Update $\bar{x}_{g0}^{i+1} \gets \bar{x}_{g,\Delta t}^i$, $\bar{p}_{g0}^{i+1} \gets \bar{p}_{g,\Delta t}^i$
\IF{$\bar{x}_{g,\Delta t}^{i} = 1$}
\STATE $\overline{UT}^{i+1}_g \gets \max(0, UT_g - \widehat{UT}^{i}_{g})$
\ELSE
\STATE  $\overline{DT}^{i+1}_g \gets \max(0, DT_g - \widehat{DT}^{i}_{g})$
\ENDIF
\RETURN $\left\{\bar{x}^{i+1}_{g0}, \overline{UT}^{i+1}_g, \overline{DT}^{i+1}_g, \bar{p}^{i+1}_{g0}\right\}$.
\end{algorithmic}
\end{algorithm}

\section{Dynamic Enforcement for security constraints}
\label{sec:dynamic_cut}

Large-scale SCUC problems face severe computational challenges due to the exponential growth of dense security constraints. Existing approaches~\cite{ardakani2013identification,xavier2019transmission} mitigate this by iteratively identifying and adding violated constraints. While these methods significantly reduce problem size by activating only critical constraints, they remain computationally expensive as they require repeatedly solving the MILP model from scratch. This sequential process delays obtaining even an initial feasible solution until complete convergence.
In this work, we propose a more effective strategy that dynamically enforces security constraints within the branch-and-cut procedure. Our key innovation involves real-time constraint handling: whenever an integer solution at a branch-and-bound node is found, we immediately verify all security constraints and activate the most violated ones without restarting the branch-and-cut process. This approach maintains the solver's internal state while selectively enforcing only the most critical constraints at each branch-and-bound node.
By integrating cut separation natively into the optimization workflow, our method could potentially achieve substantial computational savings compared to traditional iterative approaches, while preserving solution accuracy.

\begin{algorithm}[htbp]
\caption{Dynamic Cut Separation }\label{alg:subproblem}
\begin{algorithmic}[1]
\STATE \textbf{Input:} $\mathcal{P}(\emptyset)$, security constraint set $\bar{\Gamma}$.
\STATE \textbf{Initialize:} $\Gamma \gets \emptyset$, BKS = $\emptyset$.

\WHILE{Branch-and-Cut for $\mathcal{P}(\Gamma)$ not terminate}
\IF{Found integer solution $\xi$}
    \STATE $\Gamma^{\text{vio}}\gets $ Violated constraints in $\bar{\Gamma}$ at $\xi$  \hfill  \textit{ // in parallel} \label{line:viol}
    
    \IF{$\Gamma^{\text{vio}} \not = \emptyset$}
        \STATE $\Gamma\gets \Gamma\cup \Gamma^{\text{vio}}$ \hfill  \textit{// add cuts via callback}
    \ELSE
         \STATE Update BKS  \hfill  \textit{// best known solution
} 
    \ENDIF
\ENDIF
\ENDWHILE
\RETURN BKS  
\end{algorithmic}
\end{algorithm}
A pseudocode illustrating the use of branch-and-cut with dynamic cut separation to solve the subproblems is summarized in Algorithm~\ref{alg:subproblem}.
Let $\mathcal{P}(\bar{\Gamma})$ represent the SCUC model for a given planning horizon, where $\bar{\Gamma}$ represents the complete set of associated security constraints. For clarity, we omit time-interval subscripts.
The solution process begins by solving a relaxed model $\mathcal{P}(\emptyset)$ that initially omits all security constraints. Within the branch-and-cut procedure, the solver dynamically identifies violated security constraints $\Gamma^{\text{vio}} \subseteq \bar{\Gamma}$ through callback functions whenever an integer solution is found. These violated constraints are then added to the active model. If no violations are detected, the incumbent solution may be updated with the current feasible solution.
The branch-and-cut procedure with dynamic cut separation terminates upon reaching the specified time limit or achieving a predefined optimality gap. 
If a subproblem is determined to be infeasible, Algorithm~\ref{alg:subproblem} returns $\emptyset$ to signify the absence of feasible solutions.
The integration of dynamic cut separation with the solver's branch-and-cut procedure eliminates the need to iteratively solve new MILP problems from scratch, leading to greater efficiency.
To further enhance performance, the identification of these violated constraints (Step~\ref{line:viol}) can be parallelized using multi-threading.

\section{Temporal decomposition with relax-and-cut}
\label{sec:overall_algo}
In this section, we first summarize the integrated relax-and-cut framework and then introduce the RINS strategy to further refine SCUC solutions.

\subsection{Overall Algorithm}
Let $\mathcal{P}\bigl(w_I, w_R; w_F, \xi_F, \{\Gamma_t\}_{t\in w_I\cup w_R}\bigr)$ denote the optimization problem with variables and constraints associated with windows $w_I$ (Integer) and $w_R$ (Relaxed), where problem parameters are updated using fixed values $\xi_F$ from the preceding fixed window $w_F$. Here, $\{\Gamma_t\}_{t\in w_I\cup w_R}$ represents the security constraints for time intervals in $w_I\cup w_R$. This formulation is algebraically equivalent to $\mathrm{SCUC}^{PR}$, allowing subproblem updates via Algorithm~\ref{alg:update}.
The full SCUC problem corresponds to the special case $\mathcal{P}(\mathcal{T}, \emptyset; \emptyset, \emptyset, \bar{\Gamma})$, consistent with formulation~(\ref{eq:param_prob}). 

The complete relax-and-cut methodology is formalized in Algorithm~\ref{alg:framework}. The computational procedure begins by initializing the fixed ($w_F$), integer ($w_I$), and relaxed ($w_R$) windows according to the prescribed sizes $\bar{s}_I$ and $\bar{s}_R$ for the integer and relaxed windows, respectively. 
At each iteration, the framework formulates and solves a partially relaxed subproblem using the dynamic cut separation method specified in Algorithm~\ref{alg:subproblem}. The solution segment $\xi_{[1:\Delta t]}$, corresponding to the first $\Delta t$ intervals, is appended to the fixed solution vector $\xi_F$. Subsequent subproblems are generated by advancing the integer and relaxed windows while updating parameters through Algorithm~\ref{alg:update}. This iterative process continues until the relaxed window $w_R$ reaches the termination point of the planning horizon.
The procedure concludes by solving a final subproblem where $\xi_F$ remains fixed and the integer window $w_I$ encompasses all remaining uncommitted intervals, employing dynamic cut separation to ensure constraint satisfaction. Throughout execution, encountering an infeasible subproblem triggers a restart after expanding the integer window size by the increment $\Delta s$.

\begin{algorithm}[htbp]
\caption{Temporal Decomposition with Relax-and-Cut}\label{alg:framework}
\begin{algorithmic}[1]
\STATE \textbf{Input:} Problem instance, $T$, $\bar{s}_I, \bar{s}_R$, $\Delta t$ and $\Delta s$
\STATE \textbf{Initialize:} $s_F\gets -\Delta t, s_I\gets \bar{s}_I, s_R\gets\bar{s}_R, \xi_F \gets \emptyset$

\WHILE{$s_I \leq T$}
\WHILE{$s_F + s_I + s_R \leq T$}
    \STATE \textit{// Update windows}
    \STATE $s_F \gets s_F+\Delta t$
    \STATE $w_F \gets [1, s_F]$
    \STATE $w_I \gets [s_F+1, s_F+s_I]$
    \STATE $w_R \gets [s_F+s_I+1, s_F+s_I+s_R]$

    \STATE
    \STATE \textit{// Solve a subproblem}
    \STATE $\tilde{\mathcal{P}} \gets \mathcal{P}(w_I, w_R;w_F, \xi_F,\{\Gamma_t\}_{t\in w_I\cup w_R})$
    \STATE $\xi \gets $ solution of $\tilde{\mathcal{P}}$ returned by Algorithm~\ref{alg:subproblem}.
    \IF{$\xi=\emptyset$}
        \STATE ${\xi_F} = \emptyset$, break
    \ELSE
        \STATE $\xi_F \gets \text{Append}\left(\xi_F, \xi_{[1:\Delta t]}\right)$ \hfill  \textit{// Fix more commitments}
    \ENDIF
\ENDWHILE

\STATE
\STATE \textit{// Complete remaining intervals}
\IF{$\xi_F\not = \emptyset$}
\STATE $\bar{\xi} \gets$ Solve $\mathcal{P}(\mathcal{T}\setminus w_F, \emptyset; w_F, \xi_F, \bar{\Gamma})$ with Algorithm~\ref{alg:subproblem}.
\ENDIF

\STATE
\STATE \textit{// Enlarge $s_I$}
\IF{$\bar{\xi} = \emptyset$ or $\xi_F=\emptyset$}
    \STATE $s_I\gets \min(s_I + \Delta s, T)$
    \STATE $s_F \gets -\Delta t$
\ELSE
    \RETURN $\bar{\xi}$
\ENDIF
\ENDWHILE
\end{algorithmic}
\end{algorithm}

\subsection{Solution Improvement}

\begin{algorithm}[htbp]
\caption{RINS}\label{alg:rins}
\begin{algorithmic}[1]
\STATE \textbf{Input:} Problem instance, $\mathcal{T}$ and feasible solution $\bar{\xi}$.
\WHILE{termination criteria not met}
    \STATE Select an \enquote{Integer} window $w_I$ for refinement
    \STATE $w_F\gets\mathcal{T}\setminus w_I$
    \STATE $\bar{\xi} \gets $ Solve $\mathcal{P}(w_I, \emptyset; w_F, \bar{\xi}, \bar{\Gamma})$ with Algorithm~\ref{alg:subproblem}
\ENDWHILE
\RETURN $\bar{\xi}$
\end{algorithmic}
\end{algorithm}

\begin{figure}[htbp]
    \centering
    \includegraphics[width=1\linewidth]{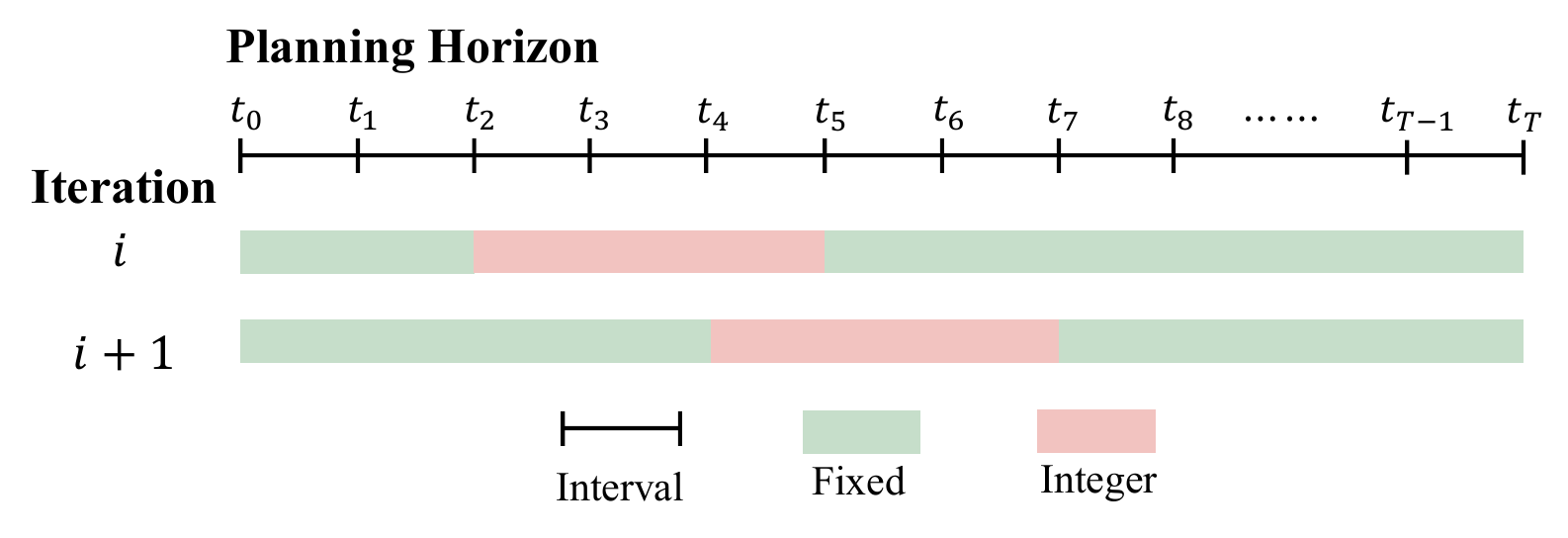}
    \caption{An illustration of iterative solution improvement.
    At the $i^{th}$ iteration, a rolling window spanning interval $(t_2, t_5]$ is dynamically selected.
    Within this window, all integer variables retain their original combinatorial freedom, while variables outside the window are fixed to their values in the incumbent solution.
    During the subsequent iteration, the window adaptively shifts to $(t_4, t_7]$.
    The process terminates upon satisfying either predefined optimality or time limits.}
    \label{fig:improvement}
\end{figure}

Given a feasible solution, we implement a time-aware refinement procedure that adapts the RINS framework for progressive improvement (see Algorithm~\ref{alg:rins}). Our method dynamically selects a moving temporal window across the planning horizon, maintaining exact MIP formulations for integer variables within the window while fixing external variables to their incumbent values. Through iterative window resizing and RINS-inspired variable unfixing, the approach systematically explores solution neighborhoods while preserving feasibility via constraint inheritance from the initial solution. 
As illustrated in Fig.~\ref{fig:improvement}, this process continues until either meeting predefined optimality thresholds or reaching computational time limits, ensuring practical deployability. The resulting refinement strategy combines global solution perspectives with localized optimization intensity, achieving quality improvements without sacrificing feasibility.

\section{Case study}
\label{sec:case_study}
This section presents a comprehensive performance evaluation of the proposed \textit{relax-and-cut temporal decomposition} method (denoted as \texttt{TD-R+Cut}) through four key analyses:
(i) Comparative benchmarking against a state-of-the-art temporal decomposition method and the monolithic SCUC model solved with Gurobi,
(ii) Ablation studies quantifying the individual contributions of the decomposition strategy and dynamic security constraint enforcement,
(iii) Parallel computing scalability assessment via thread-count sensitivity analysis, and
(iv) Performance evaluation of the proposed improvement strategy.
Throughout the experiments, \texttt{TD-R+Cut} is configured with parameters $\bar{s}_I=\bar{s}_R=\Delta t = 6$, $\Delta s = 2$.

\noindent \textbf{Baseline and Ablation Methods:}
To assess the efficacy of the proposed framework, we conduct a comparative analysis against the following baseline methods. Throughout the experiments, all problems and subproblems are solved using the state-of-the-art MILP solver Gurobi.
\begin{itemize}
    \item \texttt{Full+TF}: The \textit{full} monolithic SCUC model~(\ref{eq:param_prob}) is solved by Gurobi and the \textit{transmission constraint filtering}~(denoted as \texttt{TF}) technique from~\cite{xavier2019transmission}  is integrated to handle security constraints.
    
    \item \texttt{TD+TF}: A \textit{temporal decomposition} approach (denoted as \texttt{TD}) implemented in UnitCommitment.jl~\cite{xavier2022unitcommitment}, wherein the problem is partitioned into sequential subproblems and similarly employs \texttt{TF} from~\cite{xavier2019transmission}. 
    This method represents a degenerate case of our framework with $\bar{s}_R=0$.
    In the experiments, the window size and advancement step size are chosen to be 6. 
    
    \item \texttt{TD-E+TF}: The \texttt{TD+TF} method is used, but with a window size \textit{extended} to 12 to match the subproblem size of \texttt{TD-R+Cut}.
\end{itemize}

To further elucidate the contributions of key components within our framework, we introduce two ablation variants: 
\begin{itemize}
    \item \texttt{Full+Cut}: 
    The \textit{full} SCUC problem~(\ref{eq:param_prob}) is solved using Gurobi, and security constraints are enforced through \textit{dynamic cut separation} within the branch-and-cut procedure.

    \item \texttt{TD-R+TF}: 
    The SCUC problem is addressed using the proposed \textit{partial-relaxed temporal decomposition} method (denoted as \texttt{TD-R}), while security constraints are enforced via the \texttt{TF} method.
\end{itemize}

\noindent\textbf{Configuration.}
The proposed framework was developed as an extension of UnitCommitment.jl~\cite{xavier2022unitcommitment}.
All computational experiments were conducted using Julia v1.9.4, interfaced with the Gurobi Optimizer v11.0.2 as the underlying mathematical programming solver. 
The solver was configured with a 24-hour time limit per instance.
The MIP gap tolerance is set to $1\%$ for subproblems in \texttt{TD-R+Cut} and \texttt{TD+TF}, and to $0.1\%$ for \texttt{Full+TF} and RINS solution refinement.
The simulations were executed on a computing server featuring a 12th Gen Intel(R) Core(TM) i9-12900K processor.
All computational experiments were conducted in single-threaded mode unless otherwise noted, ensuring consistent benchmarking conditions across all tested methods.

\noindent \textbf{Benchmark.}
The proposed framework is evaluated using instances from the MATPOWER dataset~\cite{zimmerman2010matpower}. 
Specifically, we examine the 20 largest instances in the dataset, each containing over 1,000 buses; their detailed specifications are summarized in Table~\ref{tab:instance_info}. 
These instances range from moderately sized networks—with a few thousand buses, hundreds of generators, and thousands of transmission lines—to large-scale systems exceeding 13,000 buses, thousands of generators, and tens of thousands of lines. 
The number of contingencies scales with system size, starting at 1,288 for smaller cases. For the largest systems, such as 9241pegase and 13659pegase, the contingency count rises to 13,932, underscoring the computational complexity of these instances.

\begin{table}[ht]
\centering
\caption{Characteristics of selected instances}
\label{tab:instance_info}
\resizebox{0.9\linewidth}{!}{
\begin{tabular}{lllll}
\toprule
Instance & Buses & Generators & Lines & Contingencies \\
\midrule
1354pegase & 1,354 & 260 & 1,991 & 1,288 \\
1888rte & 1,888 & 296 & 2,531 & 1,484 \\
1951rte & 1,951 & 390 & 2,596 & 1,497 \\
2383wp & 2,383 & 323 & 2,896 & 2,240 \\
2736sp & 2,736 & 289 & 3,504 & 3,159 \\
2737sop & 2,737 & 267 & 3,506 & 3,161 \\
2746wop & 2,746 & 443 & 3,514 & 3,155 \\
2746wp & 2,746 & 457 & 3,514 & 3,156 \\
2848rte & 2,848 & 544 & 3,776 & 2,242 \\
2868rte & 2,868 & 596 & 3,808 & 2,260 \\
2869pegase & 2,869 & 510 & 4,582 & 3,579 \\
3012wp & 3,012 & 496 & 3,572 & 2,854 \\
3120sp & 3,120 & 483 & 3,693 & 2,950 \\
3375wp & 3,374 & 590 & 4,161 & 3,245 \\
6468rte & 6,468 & 1,262 & 9,000 & 6,094 \\
6470rte & 6,470 & 1,306 & 9,005 & 6,085 \\
6495rte & 6,495 & 1,352 & 9,019 & 6,060 \\
6515rte & 6,515 & 1,368 & 9,037 & 6,063 \\
9241pegase & 9,241 & 1,445 & 16,049 & 13,932 \\
13659pegase & 13,659 & 4,092 & 20,467 & 13,932 \\
\bottomrule
\end{tabular}
}
\end{table}

\subsection{Comparison Against Baselines}
\begin{table}[h!] 
    \centering
    \caption{Computational results on 20 selected instances.}
    \label{tab:res_1}
\renewcommand{\arraystretch}{1.1} 
    {\LARGE
    \resizebox{\linewidth}{!}{
        \begin{tabular}{l ll ll ll}
        \toprule
        \multirow{2}{*}{Instance}& \multicolumn{2}{c}{\texttt{Full+TF}} & \multicolumn{2}{c}{\texttt{TD+TF}}  & \multicolumn{2}{c}{\texttt{TD-R+Cut}} \\
        \cmidrule(l){2-3}\cmidrule(l){4-5}\cmidrule(l){6-7}
        & Obj. & t (sec) & Obj. (Gap) & t (sec) & Obj. (Gap)& t (sec)\\
        \midrule
1354pegase & 1.575 & 25 & 1.616 (2.57\%) & 16 & 1.587 (\textbf{0.75\%}) & \textbf{12} \\
1888rte & 2.345 & 149 & 2.391 (1.94\%) & 33 & 2.357 (\textbf{0.51\%}) & \textbf{30} \\
1951rte & 2.494 & 227 & 2.530 (1.42\%) & 22 & 2.511 (\textbf{0.68\%}) & \textbf{16} \\
2383wp & 1.369 & 17 & 1.380 (0.84\%) & 23 & 1.374 (\textbf{0.39\%}) & \textbf{11} \\
2736sp & 0.962 & 21 & 0.994 (3.33\%) & 24 & 0.972 (\textbf{1.03\%}) & \textbf{10} \\
2737sop & 0.847 & 17 & 0.893 (5.48\%) & 22 & 0.859 (\textbf{1.41\%}) & \textbf{8} \\
2746wop & 0.853 & 127 & - & - & 0.871 (\textbf{2.11\%}) & \textbf{17} \\
2746wp & 0.965 & 22 & 1.002 (3.85\%) & 22 & 0.984 (\textbf{1.97\%}) & \textbf{11} \\
2848rte & 2.427 & 547 & 2.463 (1.48\%) & 73 & 2.438 (\textbf{0.46\%}) & \textbf{50} \\
2868rte & 2.497 & 608 & 2.530 (1.31\%) & 45 & 2.504 (\textbf{0.29\%}) & \textbf{26} \\
2869pegase & 3.895 & 552 & 4.030 (3.47\%) & 101 & 3.936 (\textbf{1.05\%}) & \textbf{63} \\
3012wp & 1.204 & 55 & 1.243 (3.20\%) & 28 & 1.209 (\textbf{0.41\%}) & \textbf{19} \\
3120sp & 1.177 & 135 & 1.212 (2.95\%) & 44 & 1.191 (\textbf{1.20\%}) & \textbf{27} \\
3375wp & 1.418 & 2,447 & 1.555 (9.69\%) & 85 & 1.444 (\textbf{1.88\%}) & \textbf{51} \\
6468rte & 5.706 & 2,452 & 5.832 (2.21\%) & 462 & 5.757 (\textbf{0.90\%}) & \textbf{262} \\
6470rte & 6.640 & 26,985 & 6.766 (1.89\%) & \textbf{1,095} & 6.678 (\textbf{0.57\%}) & 2,593 \\
6495rte & 5.405 & 1,519 & 5.499 (1.73\%) & 350 & 5.431 (\textbf{0.48\%}) & \textbf{265} \\
6515rte & 5.573 & 1,870 & 5.689 (2.09\%) & 415 & 5.620 (\textbf{0.84\%}) & \textbf{337} \\
9241pegase & - & - & 11.646  & 3,292 & \textbf{11.381}  & \textbf{3,175} \\
13659pegase & 26.559 & 3,762 & 27.006 (1.68\%) & 1,008 & 26.714 (\textbf{0.58\%}) & \textbf{730} \\
\midrule
\bf Avg.  &   &281  & \quad\quad\quad (2.8\%)& 93 & \quad\quad\quad(0.9\%) & 59\\
\bottomrule
        \end{tabular}
        }
        }
\renewcommand{\arraystretch}{1} 
\end{table}

We first evaluated the performance of \texttt{TD-R+Cut} against both \texttt{Full+TF} and \texttt{TD+TF}. The results for 20 selected instances are presented in Table~\ref{tab:res_1}, which reports the total solve time (column \enquote{t (sec)}) in seconds and the final objective value (column \enquote{Obj.}
) scaled by a factor of $10^7$. The primal gaps (relative to \texttt{Full+TF}) for both \texttt{TD+TF} and \texttt{TD-R+Cut} appear in parentheses. 
Superior performance metrics-specifically, better objective values and less computation times-are highlighted in boldface for these two methods. The table's final row presents the geometric means of both solution times and average primal gaps across all test cases.

The \texttt{Full+TF} method represents a direct, monolithic optimization approach where the entire problem is passed to Gurobi. This strategy consistently yields optimal values for all benchmark instances except 9241pegase.
However, this solution quality comes at a significant computational cost. The solve times for \texttt{Full+TF} are prohibitively high and scale poorly with problem size. For instance, the time required increases from a modest 25 seconds for the 1354pegase instance to 2,452 seconds for 6468rte and a staggering 26,985 seconds (approximately 7.5 hours) for 6470rte.
Notably, the monolithic approach fails to produce any feasible solution for instance 9241pegase, as indicated by the \enquote{-} entries in the table. This limitation reflects that the \texttt{Full+TF} approach encounters fundamental scalability barriers. The overwhelming computational resources required for problems of this magnitude make the monolithic strategy untenable for large-scale SCUC problems. 

In contrast to the monolithic approach, \texttt{TD+TF} demonstrates a drastic reduction in computation time. On average, the solution time is just one-third of that required by \texttt{Full+TF}. This effect is most pronounced on the 6470rte instance, where the solution time decreases from 26,985 to 1,095 seconds—a remarkable speedup of nearly 25 times. Furthermore, \texttt{TD+TF} demonstrates its superior scalability by successfully finding a solution for the 9241pegase instance, on which \texttt{Full+TF} failed to provide any feasible solution. This computational efficiency is the principal promise of decomposition: by dividing a large problem into smaller, more tractable subproblems, the exponential growth in solution time is mitigated.
However, this speed is achieved at the cost of significant degradation in solution quality. The \enquote{Gap} column in the results reveals substantial primal gaps for \texttt{TD+TF}, ranging from 0.84\% on 2383wp to as high as 9.69\% on 3375wp. This gap primarily stems from myopic decision-making when handling inter-temporal constraints that connect decisions across the boundaries of time-based subproblems. Each subproblem is solved using only information within its local time window, leading to decisions that are locally optimal but globally suboptimal.
The practical consequence of this myopia is highlighted by the failure of \texttt{TD+TF} on the 2746wop instance, a case where other methods succeed. This specific failure underscores the critical importance of incorporating information about future periods when making decisions in the present.

The results demonstrate that the proposed \texttt{TD-R+Cut} method alleviates the key limitations observed in the baseline approaches. First, \texttt{TD-R+Cut} improves the optimality gap through its partially relaxed temporal decomposition, which adds a relaxed look-ahead window to each subproblem. The foresight provided by this larger window helps mitigate myopic decision-making and yields solutions closer to global optimality, while its relaxed nature effectively controls the subproblem's combinatorial complexity. Consequently, \texttt{TD-R+Cut} achieves a better relative gap than \texttt{TD+TF} on every instance, with primal gaps frequently falling below 1\%. Notably, on instance 3375wp, the gap is reduced dramatically from 9.69\% to just 1.88\%.

Importantly, this significant improvement in solution quality does not come at the expense of higher computational cost. In fact, \texttt{TD-R+Cut} is more efficient than \texttt{TD+TF}. By handling security constraints with dynamic cuts, the framework efficiently incorporates the look-ahead window with minimal overhead. 
In case 3375wp, for example, \texttt{TD-R+Cut} requires only 60\% of the solution time of \texttt{TD+TF} and achieves a nearly 50-fold speedup compared to \texttt{Full+TF}.
The average performance reported in Table~\ref{tab:res_1} shows that \texttt{TD-R+Cut} requires only one-fifth of the solution time of \texttt{Full+TF}, while also needing just half the solution time and achieving one-third the relative gap of \texttt{TD+TF}.
Overall, the results show that \texttt{TD-R+Cut} successfully identifies high-quality solutions in significantly less time than the alternatives.

\begin{figure}
    \centering
    \includegraphics[width=0.95\linewidth]{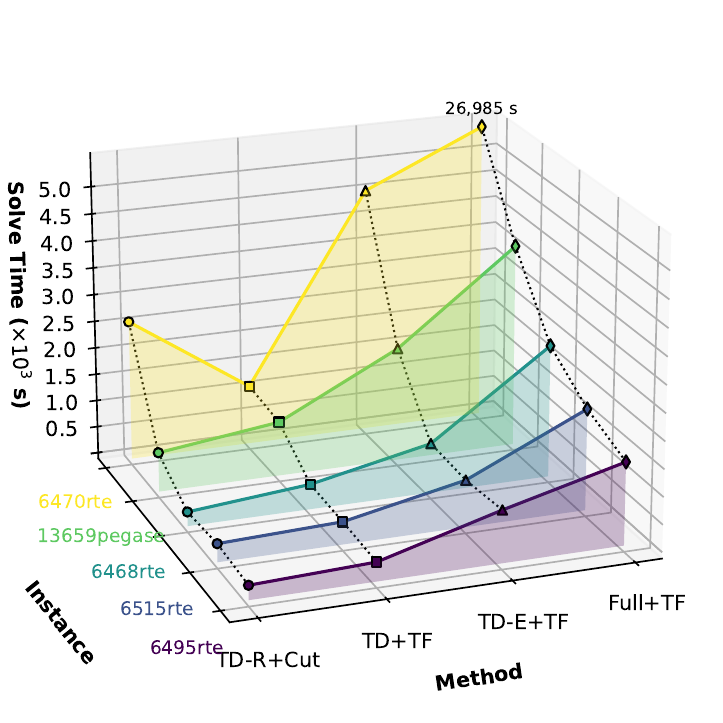}
    \caption{Comparison of solve times across five large instances. The instances are ordered by the solve time of \texttt{Full+TF}.}
    \label{fig:large_time}
\end{figure}
\begin{figure}
    \centering
    \includegraphics[width=0.95\linewidth]{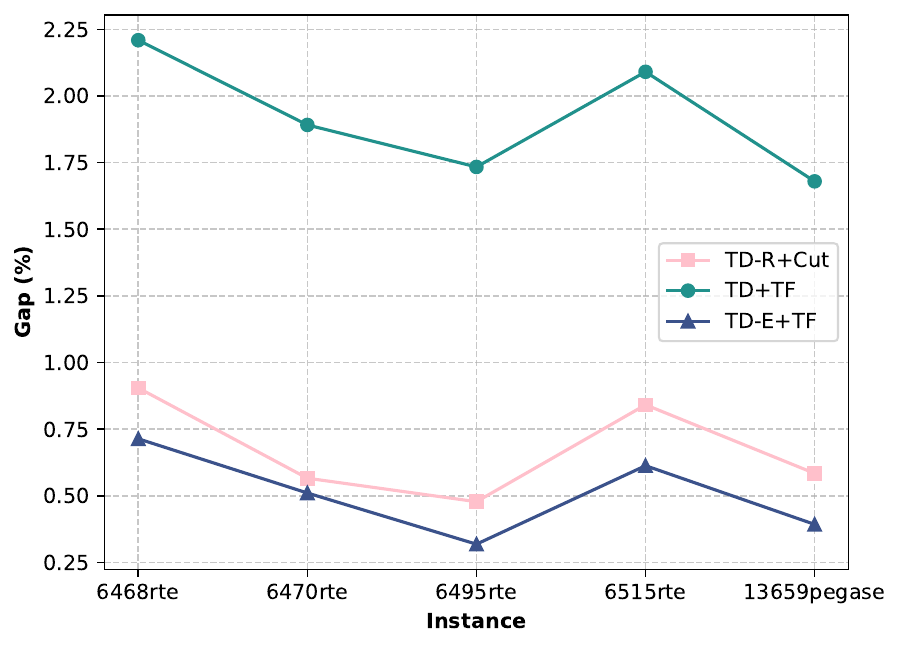}
    \caption{Comparison of relative primal gap across five large instances.}
    \label{fig:large_gap}
\end{figure}

To isolate the specific advantage of our proposed partial-relaxation mechanism, we compare \texttt{TD-R+Cut} directly against \texttt{TD-E+TF}. This alternative uses an identical subproblem size but incorporates exact lookahead information. The benefits of our approach are most apparent in the solve times for the five largest instances, plotted in Fig.~\ref{fig:large_time}. Despite using the same window size, \texttt{TD-R+Cut} achieves substantial time savings over \texttt{TD-E+TF}, as exemplified by the 13659pegase instance (730 seconds versus 2,107 seconds). This efficiency stems directly from the strategic relaxation of integer variables with less temporal dependency, which effectively expedites the MIP search process.
Crucially, this significant gain in computational speed is achieved with a negligible impact on solution quality. As shown in Fig.~\ref{fig:large_gap}, which presents the relative primal gaps, both methods consistently produce high-quality solutions with gaps well below 1.0\%. The minimal difference in solution quality between \texttt{TD-E+TF} and \texttt{TD-R+Cut} underscores the success of our partial-relaxation strategy in balancing computational tractability and solution quality.
These findings confirm that the architecture of \texttt{TD-R+Cut} is particularly advantageous for tackling large-scale SCUC problems, offering a clear path to identifying high-quality feasible solutions in a fraction of the time required by methods using exact subproblems.

\subsection{Ablation Studies}

\begin{figure*}[htp]
    \centering
    \includegraphics[width=1\linewidth]{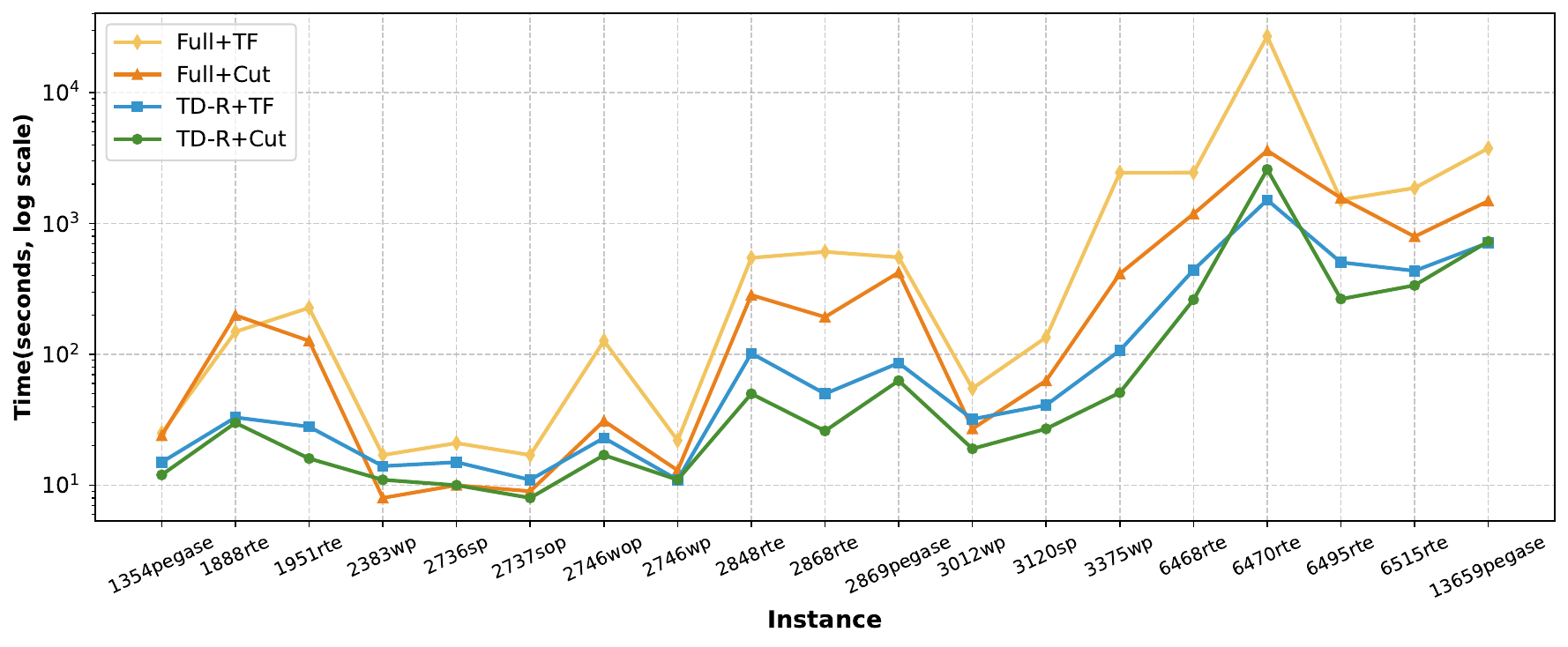}
    \caption{The ablation study compares solution times for the two techniques within \texttt{TD-R+Cut}: (i) the proposed partially relaxed temporal decomposition framework and (ii) dynamic separation for security constraints.
    The y-axis is on log-scale for better visualization.}
    \label{fig:ablation_1}
\end{figure*}

This section ablates the two key components of our proposed \texttt{TD-R+Cut} framework: partially relaxed temporal decomposition strategy and dynamic cut generation.
We compare four methods in Fig.~\ref{fig:ablation_1}: \texttt{Full+TF}, \texttt{Full+Cut}, \texttt{TD-R+TF}, and  \texttt{TD-R+Cut}, with performance measured by solution time.

The results reveal distinct performance characteristics between methods employing decomposition techniques (\texttt{TD-R+Cut}, \texttt{TD-R+TF}) and those (\texttt{Full+Cut}, \texttt{Full+TF}) solving the full optimization problem directly via Gurobi. 
While the direct approaches demonstrate competitive or even superior performance on smaller or less complex instances (e.g., 2383wp and 2736sp), likely due to solver efficiency on manageable problems, decomposition methods generally prevail.
As instance complexity increases, decomposition-based methods consistently outperform \texttt{Full+Cut} and \texttt{Full+TF}, often by a significant margin. 
The instance 6470rte presents a case where the performance difference narrows, though the trend generally holds.

The results also suggest a noticeable performance gap based on the strategy for managing security constraints. 
Methods enforcing security constraints dynamically within the branch-and-cut process, namely \texttt{TD-R+Cut} and \texttt{Full+Cut}, generally outperform their respective counterparts, \texttt{TD-R+TF} and \texttt{Full+TF}. 
This observation underscores the computational advantage of swiftly integrating the security constraints within the optimization process, thereby reducing the need for re-solving from scratch when violations are encountered.

The ablation study confirms that the partial-relaxed temporal decomposition strategy, which mitigates combinatorial complexity and the extended horizon, and the dynamic cut separation, which manages numerous security constraints, are both crucial to the effectiveness of the proposed \texttt{TD-R+Cut}.

\subsection{Parallel Performance Evaluation}
We evaluate all methods in a multi-threaded environment (8 threads maximum) to assess their practical applicability. Parallel acceleration occurs through two mechanisms: concurrent screening of violated security constraints and utilization of Gurobi's native parallel optimization capabilities. 
Table~\ref{tab:res_4} summarizes results for the six largest test cases. While all methods maintain their relative performance trends from the single-threaded case, \texttt{TD-R+Cut} demonstrates significantly better parallel scaling. On the 6515rte instance, \texttt{Full+TF} achieves a 26\% time reduction, \texttt{TD+TF} shows 14\% improvement, and \texttt{TD-R+Cut} attains 48\% acceleration - 2.8--3.4$\times$ greater gains than competing methods.
This superior parallel efficiency underscores \texttt{TD-R+Cut}'s architectural advantages for parallel computation, particularly in large-scale SCUC scenarios where rapid solution times are critical.

\begin{table}[h!]
\centering
\caption{Computational results on six large instances using 8 threads.}
\label{tab:res_4}
\resizebox{\linewidth}{!}{
{\LARGE
\renewcommand{\arraystretch}{1.1} 
    \begin{tabular}{l cl cl cl}
    \toprule
    \multirow{2}{*}{Instance}& \multicolumn{2}{c}{\texttt{Full+TF}} & \multicolumn{2}{c}{\texttt{TD+TF}}  & \multicolumn{2}{c}{\texttt{TD-R+Cut}} \\
    \cmidrule(l){2-3}\cmidrule(l){4-5}\cmidrule(l){6-7}
    & Obj. & t (sec) & Obj. (Gap) & t (sec) & Obj. (Gap)& t (sec) \\
    \midrule
6468rte & 5.706 & 793 & 5.832 (2.22\%) & 440 & 5.733 (\textbf{0.49\%}) & \textbf{143} \\
6470rte & 6.641 & 12,221 & 6.745 (1.57\%) & 777 & 6.670 (\textbf{0.44\%}) & \textbf{341} \\
6495rte & 5.403 & 2,845 & 5.494 (1.69\%) & 302 & 5.435 (\textbf{0.60\%}) & \textbf{135} \\
6515rte & 5.573 & 1,380 & 5.691 (2.12\%) & 356 & 5.613 (\textbf{0.72\%}) & \textbf{174} \\
9241pegase & - & - & \hspace{-1cm} 11.603 & \textbf{2,486} & \hspace{-1cm} \textbf{11.407} & 2,716 \\
13659pegase & 26.563 & 635 & 27.007 (1.67\%) & 1,022 & 26.734 (\textbf{0.64\%}) & \textbf{317} \\
\midrule
\bf Avg.  &   &1,891 & \quad\quad\; (1.8\%) & 674 & \quad \quad\; (0.5\%) & 315\\
\bottomrule
\end{tabular}
}
}
\renewcommand{\arraystretch}{1} 
\end{table}

\subsection{Solution Improvement}

\begin{figure}[htp]
    \centering
    \includegraphics[width=1.0\linewidth]{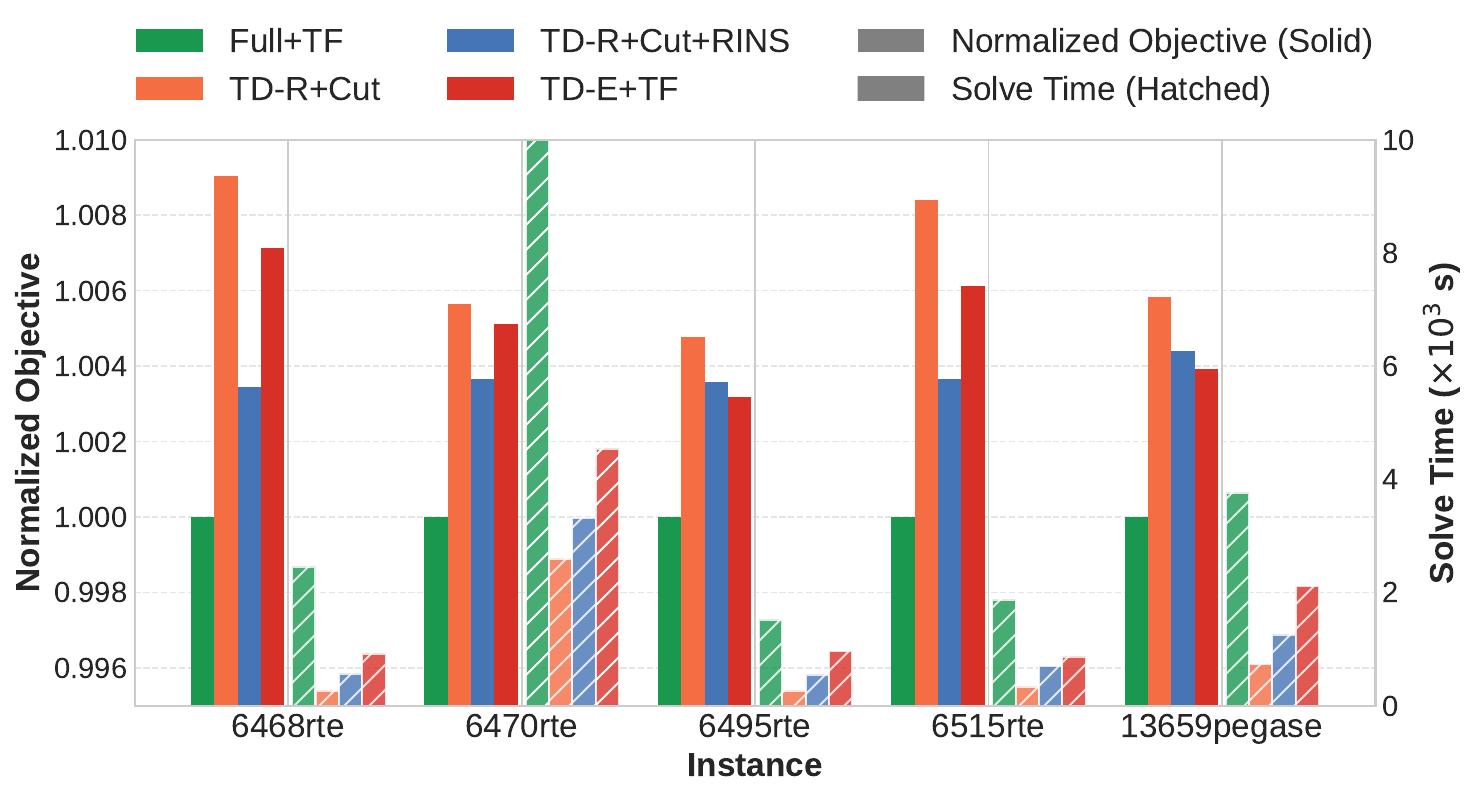}
    \caption{Comparison of objective value and solve time of \texttt{TD-R+Cut} and \texttt{TD-R+Cut+RINS} against baseline methods on five large instances.
    }
    \label{fig:improv}
\end{figure}

The \texttt{TD-R+Cut} method demonstrates strong performance while maintaining computational efficiency that enables additional solution refinement. We apply a RINS-inspired heuristic to the initial feasible solution using a sliding window strategy (12-interval window advanced by 9 intervals per iteration).
 Fig.~\ref{fig:improv} compares \texttt{Full+TF}, \texttt{TD-R+Cut}, the refined \texttt{TD-R+Cut+RINS}, and \texttt{TD-E+TF}, where solid bars indicate normalized objective values (relative to \texttt{Full+TF}) and hatched bars denote total solve times. 
Comparative analysis reveals three key advantages of our approach. First, \texttt{TD-R+Cut} requires only 15\% of the solve time needed by the monolithic \texttt{Full+TF} approach. Second, incorporating RINS refinement reduces the average primal gap from 0.67\% to 0.37\%, while total solve time remains at just 28\% of \texttt{Full+TF}'s duration. Third, the complete \texttt{TD-R+Cut+RINS} method outperforms \texttt{TD-E+TF}, delivering better solution quality (0.37\% versus 0.51\% gap) with superior computational efficiency.
These results validate our two-phase methodology combining initial solution generation via the relax-and-cut framework with targeted RINS-based refinement. The approach achieves an excellent trade-off, where a modest 13\% runtime increase yields a 45\% reduction in optimality gap (from 0.67\% to 0.37\%), making it particularly suitable for large-scale SCUC implementations that demand both solution quality and computational efficiency.

\section{conclusions}
\label{sec:conclusion}
This paper presents a novel relax-and-cut temporal decomposition strategy for efficiently obtaining high-quality SCUC solutions in large-scale power systems. 
The method successfully strikes a balance between solution quality and subproblem complexity by categorizing decision variables into \enquote{Fixed}, \enquote{Integer}, \enquote{Relaxed}, and \enquote{Ignored} across progressing time windows.
The dense security constraints are added on-the-fly within a branch-and-cut algorithm, thereby avoiding explicit enumeration.
Building upon the efficiently obtained initial solution, a RINS-inspired strategy can be subsequently employed to further refine the incumbent, if time permits. 
Numerical experiments demonstrate that the framework significantly outperforms Gurobi and state-of-the-art decomposition algorithms in terms of computational efficiency and solution quality. The results highlight its potential for real-world applications, particularly in power systems with stringent security requirements.
A current limitation of our approach is that the underlying SCUC subproblems still rely on general-purpose solvers without problem-specific customization. Developing specialized solution techniques tailored to the unique structure of these subproblems presents a promising research direction.

\bibliographystyle{IEEEtran}
\bibliography{ref}

\end{document}